\documentclass{tsp-short}
\usepackage{graphicx}

\theoremstyle{definition}

\theoremstyle{remark}

\numberwithin{equation}{section}

\begin{document}

\title
[ Estimation of  the parameters of the  Ornstein-Uhlenbeck's stochastic process]
{Estimation of  the parameters of the  Ornstein-Uhlenbeck's stochastic process }

\author{L. D. Labadze}
\address{Department of  Mathematics, Georgian Technical University  Tbilisi DC 0175, Georgia}
\curraddr{Department of  Mathematics, Georgian Technical University  Tbilisi DC 0175, Georgia}
\email{levanlabadze@yahoo.com}
\thanks{}

\author{G. R. Pantsulaia}
\address{Department of  Mathematics, Georgian Technical University  Tbilisi DC 0175, Georgia}
\curraddr{Department of  Mathematics, Georgian Technical University  Tbilisi DC 0175, Georgia}
\email{gogipantsulaia@yahoo.com}
\thanks{}


\dedicatory{Dedecated to the memory  of the
Tbilisi  State University  Professor   Grigol Sokhadze}

\begin{abstract}
It is considered  Ornstein-Uhlenbeck process $ x_t  = x_0 e^{-\theta t} + \mu (1-e^{-\theta t}) + \sigma \int_0^t e^{-\theta (t-s)} dW_s$, where
$x_0 \in R$, $\theta>0$, $ \mu \in R$ and $\sigma > 0$ are parameters. By use values  $(z_k)_{k \in N}$ of  corresponding  trajectories  at a fixed  positive  moment  $t$,
a consistent estimate  of  each  unknown parameter of  the Ornstein-Uhlenbeck's  stochastic  process is constructed under assumption  that  all another  parameters are known.
\end{abstract}

\maketitle

\section{Introduction}

In mathematics, the Ornstein-Uhlenbeck process (named after Leonard Ornstein and George Eugene Uhlenbeck  joint  celebrated   work  \cite{OrnsteinUhlenbeck1930}), is a Gauss - Markov  stochastic process ( see, for example, \cite{Feller1971}, \cite{Shiryaev80}  )  that describes the velocity of a massive Brownian particle under the influence of friction.  Over time, this process tends to drift towards its long-term mean: such a process is called mean-reverting( in this context, see, for example,  \cite{Smith2010}, \cite{Yu2009} ) .

The process can be considered to be a modification of the random walk in continuous time, or Wiener process, in which the properties of the process have been changed so that there is a tendency of the walk to move back towards a central location, with a greater attraction when the process is further away from the center. The Ornstein-Uhlenbeck process can also be considered as the continuous-time analogue of the discrete-time  process.

In recent years, however, the Ornstein-Uhlenbeck process has appeared in finance as a model
of the volatility of the underlying asset price process(see, for example, \cite{Phillips2005, Vasicek1977} ).

Note that the Ornstein-Uhlenbeck process, $x_t$  satisfies the following stochastic differential equation:
$$
dx_t = \theta (\mu-x_t) dt + \sigma dW_t \eqno(1.1)
$$
where $\theta>0$, $ \mu \in R$ and $\sigma > 0$  are parameters and $W_t$ denotes the Wiener process.

The solution  of the stochastic differential equation  (1.1)   has the following  form
$$ x_t  = x_0 e^{-\theta t} + \mu (1-e^{-\theta t}) + \sigma \int_0^t e^{-\theta (t-s)} dW_s,   \eqno (2.2)$$
where $x_0$ is assumed to be constant.

The parameters  in (2.2) have the following  sense:

(i) $\mu$ represents the equilibrium or mean value supported by fundamentals (in other words, the central location);

(ii) $\sigma$  is  the degree of volatility around it caused by shocks;

(iii) $\theta$ is  the rate by which these shocks dissipate and the variable reverts towards the mean;

(iv) $x_0$ is the underlying asset price  at moment $t=0$ ( the underlying asset initial price );

(v)  $x_t$ is the underlying asset price at moment $t>0$;

There are various scientific  papers devoted to estimate of parameter  $\mu, \sigma$  and $\theta$(see, for example \cite{Yu2009}, \cite{Smith2010}). There  least-square minimization and maximum likelihood estimation
techniques are used for the  estimating   parameters $\sigma$  and  $\mu$  which  work  successfully. The same   we can not  say concerning  the estimating  the parameter $\theta$ (see, for example, \cite{Yu2009}).

The purpose of the present paper is to introduce a new  approach which by  use  values  $(z_k)_{k \in N}$ of  corresponding  trajectories  at a fixed  positive  moment  $t$,  will allows us
to construct a consistent estimate  for  each  unknown parameter of  the Ornstein-Uhlenbeck's  stochastic  process  under an assumption  that  all another  parameters are known.

The rest of the present paper is the following:

In Section 2 we consider some auxiliary notions and  facts from the theory of stochastic differential equations   and  mathematical statistics.

In Section 3  we present the constructions of  consistent estimates  for   unknown parameters   of  the Ornstein-Uhlenbeck's  stochastic  process.

In Section 4  we present simulation  of the Ornstein–Uhlenbeck's  stochastic process and some computation results.

\section{Some auxiliary  facts from the theory of stochastic differential equations  and  mathematical statistics}

\subsection{Some auxiliary  facts from the mathematical statistics}We begin this subsection by the following definition.

\medskip

\noindent{\bf Definition 2.1.1} ( \cite{KuNi74})
~ A sequence $(x_k)_{k \in N}$ of real numbers from the interval
$(a, b)$ is said to be equidistributed or uniformly distributed on
an interval $(a, b)$ if for any subinterval $[c, d]$ of  $(a, b)$
we have
$$ \lim_{n \to \infty}
n^{-1}\#(\{x_1, x_2, \cdots, x_n\} \cap [c,d])=(b-a)^{-1}(d-c),
$$ where $\#$ denotes a counting measure.

\medskip

\noindent{\bf Definition 2.1.2}  Let $\mu$ be a
probability Borel measure on $R$ and $F$ be it's corresponding
distribution function. A sequence $(x_k)_{k \in N}$ of elements of
$R$ is said to be $\mu$-equidistributed or $\mu$-uniformly
distributed on $R$ if for every interval $[a, b] (-\infty \le a <
b \le +\infty)$ we have
$$
\lim_{n \to \infty} n^{-1} \# ([a,b] \cap \{x_1, \cdots, x_n \}
)=F(b)-F(a).
$$
\medskip

\noindent{\bf Lemma 2.1.1} {\it  Let $(x_k)_{k \in N}$ be
$\ell_1$-equidistributed sequence on $(0,1)$, $F$ be a strictly
increasing  continuous distribution function on $R$ and $p$ be a
Borel probability measure on $R$ defined by $F$. Then
$(F^{-1}(x_k))_{k \in N}$ is $p$-equidistributed on $R$.
}

\begin{proof} We have
$$
\lim_{n \to \infty}n^{-1}\# ([a,b] \cap \{F^{-1}(x_1), \cdots,
F^{-1}(x_n) \} )=
$$
$$
\lim_{n \to \infty} n^{-1}\# ([F(a),F(b)] \cap \{x_1, \cdots, x_n
\})=F(b)-F(a).\eqno
$$
\end{proof}

\medskip

\noindent{\bf Corollary 2.1.1} {\it  Let $F$ be a strictly
increasing continuous distribution  function on $R$ and $p$ be a
Borel probability measure on $R$ defined by $F$. Then for a set
$D_{F} \subset R^{N}$ of all $p$-equidistributed sequences  on $R$
we have :

(i) $D_{F}= \{ (F^{-1}(x_k))_{k \in N} : (x_k)_{k \in N} \in D
\}$;

(ii)~$p^N(D_{F})=1$.}

Let  $\{ \mu_{\theta} : \theta \in R\}$  be  a family  Borel probability measures in $R$.  By  $\mu_{{\theta}}^{N}$  we  denote
the $N$-power of the measure $\mu_{\theta}$  for $ \theta \in R$.

\medskip

\noindent{\bf Definition 2.1.3}  A Borel
measurable function  $T_n : R^n \to R ~(n \in N)$ is called a
consistent estimator of a parameter $\theta$ (in the sense of  convergence  almost
everywhere) for the family $(\mu_{\theta}^N)_{\theta
\in R}$ if the following condition
$$
\mu_{\theta}^N (\{ (x_k)_{k \in N} :~(x_k)_{k \in N} \in R^N~\&~
\lim_{n \to \infty}T_n(x_1, \cdots, x_n)=\theta \})=1
$$
holds true  for each $\theta \in R$.

\medskip

\noindent{\bf Definition 2.1.4} A Borel
measurable function  $T_n : R^n \to R ~(n \in N)$ is called a
consistent estimator of a parameter $\theta$ (in the sense of
convergence in probability) for the family
$(\mu_{\theta}^N)_{\theta \in R}$   if for every $\epsilon>0$ and
$\theta \in R$ the following condition
$$
\lim_{n \to \infty} \mu_{\theta}^N (\{ (x_k)_{k \in N} :~(x_k)_{k
\in N} \in R^N~\&~ | T_n(x_1, \cdots, x_n)-\theta|>\epsilon \})=0
$$
holds.
\medskip

\noindent{\bf Definition 2.1.5}  A Borel
measurable function  $T_n : R^n \to R ~(n \in N)$ is called a
consistent estimator of a parameter $\theta$ (in the sense of
convergence in distribution ) for the family
$(\mu_{\theta}^N)_{\theta \in R}$ if for every continuous bounded
real valued function $f$ on $R$ the following condition
$$
\lim_{n \to \infty} \int_{R^N}f(T_n(x_1, \cdots, x_n))d
\mu_{\theta}^N((x_k)_{k \in N})=f(\theta)
$$
holds.
\medskip

\noindent{\bf Remark 2.1.1}  Following
\cite{Shiryaev80} (see, Theorem 2, p. 272), for the family
$(\mu_{\theta}^N)_{\theta \in R}$ we have:

(a) an existence of a consistent estimator of a parameter $\theta$
in the sense of   convergence  almost everywhere  implies an existence of a
consistent estimator of a parameter $\theta$ in the sense of
convergence in probability;

(b) an existence of a consistent estimator of a parameter $\theta$
in the sense of convergence in probability implies an existence of
a consistent estimator of a parameter $\theta$ in the sense of
convergence in distribution.
\medskip

\noindent{\bf Definition 2.1.6}  Following
\cite{Ibram80}, the family $(\mu_{\theta}^N)_{\theta \in R}$ is
called  strictly separated if there exists a  family
$(Z_{\theta})_{\theta \in R}$ of Borel subsets of $R^N$ such  that

(i)~$\mu_{\theta}^N(Z_{\theta})=1$ for $\theta \in R$;

(ii)~$Z_{\theta_1} \cap Z_{\theta_2}=\emptyset$ for all different
parameters $\theta_1$ and $\theta_2$ from $R$.

(iii)~$\cup_{\theta \in R}Z_{\theta}=R^N.$

\medskip

\noindent{\bf Definition 2.1.7}  Following
\cite{Ibram80}, a Borel measurable function $T : R^N \to R$ is
called  an infinite sample  consistent  estimator of a parameter $\theta$ for
the family $(\mu_{\theta}^N)_{\theta \in R}$ if the following
condition
$$
(\forall \theta)(\theta \in R \rightarrow \mu_{\theta}^N (\{
(x_k)_{k \in N}~ : ~(x_k)_{k \in N} \in R^N ~\&~ T((x_k)_{k \in
N})=\theta\})=1)
$$
holds.
\medskip

\noindent{\bf Remark 2.1.2}  Note that an
existence of an infinite sample consistent estimator of a parameter $\theta$
for the family $(\mu_{\theta}^N)_{\theta \in R}$ implies that the
family $(\mu_{\theta}^N)_{\theta \in R}$ is strictly separated.
Indeed, if we  set $Z_{\theta}=\{ (x_k)_{k \in N} : (x_k)_{k \in
N} \in R^N ~\&~T((x_k)_{k \in N})=\theta\}$ for $\theta \in R$,
then all conditions in Definition 2.1.6  will be satisfied.

In the sequel we will need the well known fact from the
probability theory (see, for example, \cite{Shiryaev80}, p. 390).
\medskip

\noindent{\bf Lemma 2.1.2}  (Kolmogorov's strong law of large numbers)
{\it Let $X_1, X_2, ...$  be a  sequence of independent identically
distributed random variables  defined on the probability space
$(\Omega, \mathcal{F},P)$. If  these random variables have a
finite expectation $m$ (i.e., $E(X_1) = E(X_2) = ... = m <
\infty$), then the following condition
$$
P(\{ \omega : \lim_{n \to \infty} n^{-1}\sum_{k=1}^nX_k(\omega)=m
\})=1
$$
holds  true. }

\subsection{Some auxiliary  facts from the theory of stochastic differential equations}

 By use  approuches introduced in  \cite{Protter2004}  one can get the validity of the following assertions.

\medskip

\noindent{\bf Lemma 2.2.1}{\it ~Let's  consider  an Ornstein–Uhlenbeck process $x_t$  satisfies the following stochastic differential equation:
$$
dx_t = \theta (\mu-x_t) dt + \sigma dW_t \eqno (3.1.1)
$$
where $\theta>0$, $ \mu$ and $\sigma > 0$  are parameters and $W_t$ denotes the Wiener process.
Then the solution of this stochastic differential equation  (3.1.1)  is given by
$$ x_t  = x_0 e^{-\theta t} + \mu (1-e^{-\theta t}) + \sigma \int_0^t e^{-\theta (t-s)} dW_s,   \eqno(3.1.2)$$
where $x_0$ is assumed to be constant.}

\begin{proof}
The stochastic differential equation  (3.1.1)   is solved by  variation of parameters.  Changing variable

$$ f(x_t, t) = x_t e^{\theta t} $$

we get

$$
d f(x_t,t) =  \theta x_t e^{\theta t} dt + e^{\theta t} dx_t
 = e^{\theta t}\theta \mu  dt + \sigma e^{\theta t}  dW_t.
$$

Integrating from $0$ to $t$   we get

$$ x_t e^{\theta t} = x_0 + \int_0^t e^{\theta s}\theta\,\mu \, ds + \int_0^t \sigma\,e^{\theta s}\, dW_s $$

whereupon we see

$$ x_t  = x_0 e^{-\theta t} + \mu (1-e^{-\theta t}) + \sigma \int_0^t e^{-\theta (t-s)} dW_s. $$
\end{proof}
\medskip

\noindent{\bf Lemma 2.2.2}{\it ~  Under  conditions  of  Lemma 2.2.1, the following  equalities

(i)~ $E(x_t)=x_0 e^{-\theta t}+\mu(1-e^{-\theta t}); $

(ii)~   $\mbox{cov}(x_s,x_t)=\frac{\sigma^2}{2\theta}\left( e^{-\theta(t-s)} - e^{-\theta(t+s)} \right) ;$

(iii) ~$\mbox{var}(x_s) = \frac{\sigma^2}{2\theta}\left( 1 - e^{-2\theta s)} \right);$

hold true.}

\begin{proof}  The validity  of the item (i) is obvious.  In order to prove the validity of the items (ii)-(iii),  we can use the Ito  isometry  to calculate the covariance function  by
$$
\mbox{cov}(x_s,x_t) = E[(x_s - E[x_s])(x_t - E[x_t])]
 = E \left[ \int_0^s \sigma  e^{\theta (u-s)}\, dW_u \int_0^t \sigma  e^{\theta (v-t)}\, dW_v \right]
$$
$$
 = \sigma^2 e^{-\theta (s+t)}E \left[ \int_0^s  e^{\theta u}\, dW_u \int_0^t  e^{\theta v}\, dW_v \right]
= \frac{\sigma^2}{2\theta} \, e^{-\theta (s+t)}(e^{2\theta \min(s,t)}-1).
$$
Thus if $s<t$(so that $min(s,t)=s$), then we have
$$\mbox{cov}(x_s,x_t) = \frac{\sigma^2}{2\theta}\left( e^{-\theta(t-s)} - e^{-\theta(t+s)} \right).
$$
Similarly, if $s=t$ (so that $min(s,t)=s)$, then we have
$$\mbox{var}(x_s) = \frac{\sigma^2}{2\theta}\left( 1 - e^{-2\theta s)} \right).
$$

\end{proof}

\section{Estimation of parameters  of the  Ornstein - Uhlenbeck stochastic model }

\subsection{Estimation of the the underlying asset initial price $x_0$    in an Ornstein - Uhlenbeck stochastic model }

Let consider Ornstein- Uhlenbeck process

$$ x_t  = x_0 e^{-\theta t} + \mu (1-e^{-\theta t}) + \sigma \int_0^t e^{-\theta (t-s)} dW_s, \eqno  (3.1.1) $$

where $\theta>0$, $\mu \in R$,$\sigma>0$ and $W_s$ is Wiener process.

Here

(i) $\mu$ represents the equilibrium or mean value supported by fundamentals;

(ii) $\sigma$  is  the degree of volatility around it caused by shocks;

(iii) $\theta$ is  the rate by which these shocks dissipate and the variable reverts towards the mean;

(iv) $x_0$ is the underlying asset  price at initial  moment  $t=0$(the underlying asset initial price);

(v)  $x_t$ is the underlying asset  price  at moment $t~( t>0)$;

\medskip

\noindent{\bf Lemma 3.1.1}{\it~  For  $t>0$, $x_0 \in R$, $\theta>0$, $\mu \in R$  and $\sigma>0$,   let's  $\gamma_{(t,x_0.\theta,\mu,\sigma)}$ be a Gaussian probability  measure in  $R$ with the mean  $m_t=x_0 e^{-\theta t}+\mu(1-e^{-\theta t}) $ and  the variance $\sigma_t^2=\frac{\sigma^2}{2\theta}\left( 1 - e^{-2\theta s}\right)$. Assuming that  parameters $t$, $\theta$, $\mu$  and $\sigma$ are fixed,  denote by
$\gamma_{x_0}$ the measure $\gamma_{(t,x_0.\theta,\mu,\sigma)}$.
Let define the estimate $T_n : R^n  \to R$    by the following formula
$$
T_n((z_k)_{1 \le k \le n})=e^{\theta t}  \frac{ \sum_{k=1}^n z_k }{n} -\mu e^{\theta t} (1-e^{-\theta t}).    \eqno  (3.1.2)
$$
Then we get
$$
\gamma_{x_0}^{\infty}\{ (z_k)_{k \in N} :  (z_k)_{k \in N} \in  R^{\infty} ~\&~   \lim_{n \to \infty}T_n((z_k)_{1 \le k \le n})=x_0 \}=1,\eqno  (3.1.3)
$$
provided that $T_n$ is a consistent estimator of the underlying asset  price $x_0 \in R$ in the sense of  convergence almost  everywhere  for the family  of probability measures  $(\gamma_{x_0})_{x_0 \in R}$. }

\begin{proof} Let's consider probability space $(\Omega, \mathcal{F},P)$, where $\Omega=R^{\infty}$, $\mathcal{F}=B(R^{\infty})$, $P=\gamma_{x_0}^{\infty}$.

For $k \in N$ we consider $k$-th projection $Pr_k$   defined on  $R^{\infty}$  by
$$
Pr_k( (x_i)_{i \in N})=x_k\eqno  (3.1.4)
$$
for $(x_i)_{i \in N} \in R^{\infty}$.

It is obvious that $(Pr_k)_{k \in N}$ is  sequence of independent Gaussian random variables with   the mean  $m_t=x_0 e^{-\theta t}+\mu(1-e^{-\theta t}) $ and the variance  $\sigma_t^2=\frac{\sigma^2}{2\theta}\left( 1 - e^{-2\theta s} \right)$. By use  Kolmogorov  Strong  Law of Large numbers we get

$$\gamma_{x_0}^{\infty}\{ (z_k)_{k \in N} \in  R^{\infty} ~\&~   \lim_{n \to \infty}\frac{\sum_{k=1}^n Pr_k( (z_k)_{k \in N})}{n}=x_0 e^{-\theta t}+\mu(1-e^{-\theta t})  \}=1,  \eqno  (3.1.5)
$$
which implies
$$
\gamma_{x_0}^{\infty}\{ (z_k)_{k \in N} \in  R^{\infty} ~\&~   \lim_{n \to \infty}\left (
  e^{\theta t}\frac{\sum_{k=1}^n z_k}{n} -e^{\theta t}\mu(1-e^{-\theta t}) \right )   =x_0  \}$$
$$
=\gamma_{x_0}^{\infty}\{ (z_k)_{k \in N}  \in  R^{\infty} ~\&~   \lim_{n \to \infty}T_n((z_k)_{1 \le k \le  n})=x_0 \}=1.
$$

\end{proof}

\noindent{\bf Remark 3.1.1}  By  use  Remark 2.1.1  and   Lemma 3.1.1  we deduce that  $T_n$ is a consistent estimator of  the underlying asset  price $x_0 \in R$ in the sense of
convergence in probability for the statistical structure $(\gamma_{x_0})_{x_0 \in R}$  as well $T_n$ is a consistent estimator of the underlying asset  price  $x_0 \in R$ in the sense of
convergence in distribution for the statistical structure $(\gamma_{x_0})_{x_0 \in R}$.

\medskip

\noindent{\bf Theorem 3.1.1}{\it  ~Suppose that the family of probability measures $(\gamma_{x_0}^{\infty})_{x_0 \in R}$ and the estimators $T_n:R^n \to R(n \in N)$ come  from   Lemma 3.1.1.  Then the estimators
$T^{(0)}:R^{\infty} \to  R $  and  $T^{(1)}:R^{\infty} \to  R $ defined by
$$
T^{(0)}((z_k)_{k \in N})=\underline{\lim}_{n \to \infty}T_n((z_k)_{1\le k \le n}) \eqno(3.1.6)
$$
and
$$
T^{(1)}((z_k)_{k \in N})=\overline{\lim}_{n \to \infty}T_n((z_k)_{1\le k \le n}). \eqno(3.1.7)
$$
are  infinite-sample consistent  estimators of the underlying asset  price  $x_0$  for the family of probability measures $(\gamma_{x_0}^{\infty})_{x_0 \in R}$.  }

\begin{proof} Note that  we have
\begin{align*}
&\gamma_{x_0}^{\infty}\{ (z_k)_{k \in N}  \in  R^{\infty} ~\&~  T^{(0)}((z_k)_{k \in N})=x_0 \}\\
&=\gamma_{x_0}^{\infty}\{ (z_k)_{k \in N}  \in  R^{\infty} ~\&~ \underline{\lim}_{n \to \infty}T_n((z_k)_{1 \le k \le n})=x_0 \}\\
& \ge \gamma_{x_0}^{\infty}\{ (z_k)_{k \in N}  \in  R^{\infty} ~\&~   \lim_{n \to \infty}T_n((z_k)_{1 \le k \le n})=x_0 \}=1,\\
\end{align*}
which means  that  $T^{(0)}$ is an infinite-sample consistent  estimator of the underlying asset  price  $x_0$  for the family of probability measures $(\mu_{x_0}^{\infty})_{x_0 \in R}$.

Similarly, we have

\begin{align*}
&\gamma_{x_0}^{\infty}\{ (z_k)_{k \in N} \in  R^{\infty} ~\&~  T^{(1)}((z_k)_{k \in N})=x_0 \}\\
&=\gamma_{x_0}^{\infty}\{ (z_k)_{k \in N} \in  R^{\infty} ~\&~ \overline{\lim}_{n \to \infty}T_n((z_k)_{1 \le k \le n})=x_0 \}\\
& \ge \gamma_{x_0}^{\infty}\{ (z_k)_{k \in N}  \in  R^{\infty} ~\&~   \lim_{n \to \infty}T_n((z_k)_{(z_k)_{1 \le k \le n}})=x_0 \}=1,
\end{align*}
which means  that  $T^{(1)}$ is an infinite-sample consistent  estimator of the underlying asset  pricel  $x_0$  for the family of probability measures $(\gamma_{x_0}^{\infty})_{x_0 \in R}$.

\end{proof}

\noindent{\bf Remark 3.1.2}  By use  Remark 2.1.2  we deduce  that an
existence of  infinite sample consistent estimators  $T^{(0)}$  and  $T^{(1)}$ of the underlying asset  price $x_0$
for the family $(\gamma_{x_0}^{\infty})_{x_0 \in R}$ (cf. Theorem 3.1.1) implies that the
family $(\gamma_{x_0}^{\infty})_{x_0 \in R}$ is strictly separated.

\subsection{Estimation of  the equilibrium  $\mu$   in an Ornstein - Uhlenbeck stochastic model }

This subsection we begin by the following assertion.
\medskip

\noindent{\bf Lemma 3.2.1}  {\it  For  $t>0$, $x_0 \in R$, $\theta>0$, $\mu \in R$  and $\sigma>0$,   let's  $\gamma_{(t,x_0.\theta,\mu,\sigma)}$ be a Gaussian probability  measure in  $R$ with the mean  $m_t=x_0 e^{-\theta t}+\mu(1-e^{-\theta t}) $ and  the variance $\sigma_t^2=\frac{\sigma^2}{2\theta}\left( 1 - e^{-2\theta s}\right)$. Assuming that  parameters $x_0$, $t$, $\theta$  and $\sigma$ are fixed,  for $\mu \in R$  let's  denote by $\gamma_{\mu}$ the measure $\gamma_{(t,x_0.\theta,\mu,\sigma)}$.
Let define the estimate $T^*_n : R^n  \to R$    by the following formula
$$
T^*_n((z_k)_{1 \le k \le n})= \left(\frac{ \sum_{k=1}^n z_k }{n} -  x_0 e^{-\theta t}\right)  /  \left (1-e^{-\theta t} \right).    \eqno  (3.2.1)
$$
Then we get
$$
\gamma_{\mu}^{\infty}\{ (z_k)_{k \in N} :  (z_k)_{k \in N} \in  R^{\infty} ~\&~   \lim_{N \to \infty}T_n((z_k)_{1 \le k \le n})=x_0 \}=1,\eqno  (3.2.2)
$$
provided that $T_n$ is a consistent estimator of the equilibrium  $\mu \in R$ in the sense of    convergence  almost  everywhere for the family of  probability  measures $(\gamma_{\mu}^{\infty})_{\mu  \in R}$. }

\begin{proof} Let's consider probability space $(\Omega, \mathcal{F},P)$, where $\Omega=R^{\infty}$, $\mathcal{F}=B(R^{\infty})$, $P=\gamma_{\mu}^{\infty}$.

For $k \in N$ we consider $k$-th projection $Pr_k$   defined on  $R^{\infty}$  by
$$
Pr_k( (x_i)_{i \in N})=x_k\eqno  (3.2.3)
$$
for $(x_i)_{i \in N} \in R^{\infty}$.

It is obvious that $(Pr_k)_{k \in N}$ is  sequence of independent Gaussian random variables  with the mean  $m_t=x_0 e^{-\theta t}+\mu(1-e^{-\theta t}) $ and the variance  $\sigma_t^2=\frac{\sigma^2}{2\theta}\left( 1 - e^{-2\theta s} \right)$. By use  Kolmogorov's  Strong  Law of Large numbers we get

$$\gamma_{\mu}^{\infty}\{ (z_k)_{k \in N} \in  R^{\infty} ~\&~   \lim_{n \to \infty}\frac{\sum_{k=1}^n Pr_k( (z_k)_{k \in N})}{n}=x_0 e^{-\theta t}+\mu(1-e^{-\theta t})  \}=1,  \eqno  (3.2.4)
$$
which implies
$$
\gamma_{\mu}^{\infty}\{ (z_k)_{k \in N} \in  R^{\infty} ~\&~   \lim_{n \to \infty}\left(\frac{ \sum_{k=1}^n z_k }{n} -  x_0 e^{-\theta t}\right)  /  \left (1-e^{-\theta t} \right)  =\mu  \}$$
$$
=\gamma_{\mu}^{\infty}\{ (z_k)_{k \in N}  \in  R^{\infty} ~\&~   \lim_{n \to \infty}T^*_n((z_k)_{1 \le k \le  n})=\mu\}=1.
$$

\end{proof}

\noindent{\bf Remark 3.2.1}  By  use  Remark 2.1.1  and   Lemma 3.2.1  we deduce that  $T^*_n$ is a consistent estimator of   the equilibrium  $\mu \in R$   in the sense of
convergence in probability for the family  of measures  $(\gamma_{\mu}^{\infty})_{\mu \in R}$  as well $T_n$ is a consistent estimator of  the equilibrium  $\mu \in R$   in the sense of
convergence in distribution for  the family  of measures  $(\gamma_{\mu}^{\infty})_{\mu \in R}$.

\medskip

\noindent{\bf Theorem 3.2.1}{\it ~Suppose that the family of probability measures $(\gamma_{\mu}^{\infty})_{\mu \in R}$ and the estimators $T^*_n:R^n \to R(n \in N)$ come  from   Lemma 3.2.1.  Then the estimators
$T_*^{(0)}:R^{\infty} \to  R $  and  $T_{*}^{(1)}:R^{\infty} \to  R $ defined by
$$
T_*^{(0)}((z_k)_{k \in N})=\underline{\lim}_{n \to \infty}T^*_n((z_k)_{1\le k \le n}) \eqno(3.2.5)
$$
and
$$
T_*^{(1)}((z_k)_{k \in N})=\overline{\lim}_{n \to \infty}T^*_n((z_k)_{1\le k \le n}). \eqno(3.2.6)
$$
are  infinite-sample consistent  estimators of the equilibrium  $\mu \in R$   for the family of probability measures  $(\gamma_{\mu}^{\infty})_{\mu \in R}$.  }

\begin{proof} Note that  we have
\begin{align*}
&\gamma_{\mu}^{\infty}\{ (z_k)_{k \in N}  \in  R^{\infty} ~\&~  T_{*}^{(0)}((z_k)_{k \in N})=\mu \}\\
&=\gamma_{\mu}^{\infty}\{ (z_k)_{k \in N}  \in  R^{\infty} ~\&~ \underline{\lim}_{n \to \infty}T^*_n((z_k)_{1 \le k \le n})=\mu \}\\
& \ge \gamma_{\mu}^{\infty}\{ (z_k)_{k \in N}  \in  R^{\infty} ~\&~   \lim_{N \to \infty}T^*_n((z_k)_{1 \le k \le n})=\mu \}=1,\\
\end{align*}
which means  that  $T_{*}^{(0)}$ is an infinite-sample consistent  estimators of  the equilibrium  $\mu \in R$   for the family of probability measures$(\gamma_{\mu}^{\infty})_{\mu \in R}$.

Similarly, we have

\begin{align*}
&\gamma_{\mu}^{\infty}\{ (z_k)_{k \in N} \in  R^{\infty} ~\&~  T_{*}^{(1)}((z_k)_{k \in N})=\mu \}\\
&=\gamma_{\mu}^{\infty}\{ (z_k)_{k \in N} \in  R^{\infty} ~\&~ \overline{\lim}_{n \to \infty}T^*_n((z_k)_{1 \le k \le n})=\mu\}\\
& \ge \gamma_{\mu}^{\infty}\{ (z_k)_{k \in N}  \in  R^{\infty} ~\&~   \lim_{N \to \infty}T^*_n((z_k)_{(z_k)_{1 \le k \le n}})=\mu \}=1,
\end{align*}
which means  that  $T_{*}^{(1)}$ is an infinite-sample consistent  estimators of  the equilibrium  $\mu \in R$   for the family of probability measures$(\gamma_{\mu}^{\infty})_{\mu \in R}$.

\end{proof}

\noindent{\bf Remark 3.2.2}  Note that an
existence of  infinite sample consistent estimators  $T_{*}^{(0)}$  and  $T_{*}^{(1)}$  of  the equilibrium  $\mu \in R$   for the family of probability measures $(\gamma_{\mu}^{\infty})_{\mu \in R}$
 (cf. Theorem 3.2.1) implies that the
family $(\gamma_{\mu}^{\infty})_{\mu \in R}$   is strictly separated.

\subsection{Estimation of  the rate  $\theta$  in  Ornstein - Uhlenbeck stochastic model }

We begin this subsection by the following  lemma.

\noindent{\bf Lemma 3.3.1}  {\it  For  $t>0$, $x_0 \in R$, $\theta>0$, $\mu \in R$  and $\sigma>0$,   let's  $\gamma_{(t,x_0.\theta,\mu,\sigma)}$ be a Gaussian probability  measure in  $R$ with the mean  $m_t=x_0 e^{-\theta t}+\mu(1-e^{-\theta t}) $ and  the variance $\sigma_t^2=\frac{\sigma^2}{2\theta}\left( 1 - e^{-2\theta s}\right)$. Assuming that  parameters $x_0$, $t$, $\mu$  and $\sigma$ are fixed  such that $x_0\neq \mu$,  for $\theta >0$,  let's  denote by $\gamma_{\theta}$ the measure $\gamma_{(t,x_0,\theta,\mu,\sigma)}$.
Let define the estimate $T^{**}_n : R^n  \to R$    by the following formula
$$
T^{**}_n((z_k)_{1 \le k \le n})= -\frac{1}{t}  \ln  \left(  \frac{\frac{\sum_{k=1}^n z_k}{n}-\mu}{  x_0 -\mu}  \right).    \eqno  (3.2.2)
$$
Then we get
$$
\gamma_{\theta}^{\infty}\{ (z_k)_{k \in N} :  (z_k)_{k \in N} \in  R^{\infty} ~\&~   \lim_{N \to \infty}T^{**}_n((z_k)_{1 \le k \le n})=\theta \}=1,\eqno  (3.2.3)
$$
provided that $T^{**}_n$ is a consistent estimator of the rate $\theta>0$ in the sense of  convergence   almost  everywhere   for the family of  probability  measures $(\gamma_{\theta}^{\infty})_{\theta >0}$. }

\begin{proof} Let's consider probability space $(\Omega, \mathcal{F},P)$, where $\Omega=R^{\infty}$, $\mathcal{F}=B(R^{\infty})$, $P=\gamma_{\theta}^{\infty}$.

For $k \in N$ we consider $k$-th projection $Pr_k$   defined on  $R^{\infty}$  by
$$
Pr_k( (x_i)_{i \in N})=x_k\eqno  (3.2.4)
$$
for $(x_i)_{i \in N} \in R^{\infty}$.

It is obvious that $(Pr_k)_{k \in N}$ is  sequence of independent Gaussian random variables  with the mean  $m_t=x_0 e^{-\theta t}+\mu(1-e^{-\theta t}) $ and the variance  $\sigma_t^2=\frac{\sigma^2}{2\theta}\left( 1 - e^{-2\theta s} \right)$. By use  Kolmogorov's  Strong  Law of Large numbers  we get

$$\gamma_{\theta}^{\infty}\{ (z_k)_{k \in N} \in  R^{\infty} ~\&~   \lim_{n \to \infty}\frac{\sum_{k=1}^n Pr_k( (z_k)_{k \in N})}{n}=x_0 e^{-\theta t}+\mu(1-e^{-\theta t})  \}=1,  \eqno  (3.2.5)
$$
which implies
$$
\gamma_{\theta}^{\infty} \{ (z_k)_{k \in N} \in  R^{\infty} ~\&~   \lim_{n \to \infty}  -\frac{1}{t}  \ln  \left(  \frac{\frac{\sum_{k=1}^n z_k}{n}-\mu}{  x_0 -\mu} \right)=\theta
  \}
$$
$$
=\gamma_{\theta}^{\infty}\{ (z_k)_{k \in N}  \in  R^{\infty} ~\&~   \lim_{n \to \infty}T^{**}_n((z_k)_{1 \le k \le  n})=\theta\}=1.
$$

\end{proof}

\noindent{\bf Remark 3.3.1}  By  use  Remark 2.1.1  and   Lemma 3.3.1  we deduce that  $T^{**}_n$ is a consistent estimator of   the equilibrium  $\theta $   in the sense of
convergence in probability for the family  of  probability  measures  $(\gamma_{\theta}^{\infty})_{\theta>0}$  as well $T^{**}_n$ is a consistent estimator of  $\theta$   in the sense of
convergence in distribution for  the same family  of probability  measures.

\medskip

\noindent{\bf Theorem 3.3.1}{\it ~ Suppose that the family of probability measures $(\gamma_{\theta}^{\infty})_{\theta>0}$ and the estimators $T^{**}_n:R^n \to R(n \in N)$ come  from   Lemma 3.3.1.  Then the estimators
$T_{**}^{(0)}:R^{\infty} \to  R $  and  $T_{**}^{(1)}:R^{\infty} \to  R $ defined by
$$
T_{**}^{(0)}((z_k)_{k \in N})=\underline{\lim}_{n \to \infty}T^{**}_n((z_k)_{1\le k \le n}) \eqno(3.2.6)
$$
and
$$
T_{**}^{(1)}((z_k)_{k \in N})=\overline{\lim}_{n \to \infty}T^{**}_n((z_k)_{1\le k \le n}). \eqno(3.2.7)
$$
are  infinite-sample consistent  estimators of the rate $\theta$  for the family of probability measures  $(\gamma_{\theta}^{\infty})_{\theta>0}$.  }

\begin{proof} Note that  we have
\begin{align*}
&\gamma_{\theta}^{\infty}\{ (z_k)_{k \in N}  \in  R^{\infty} ~\&~  T_{**}^{(0)}((z_k)_{k \in N})=\mu \}\\
&=\gamma_{\theta}^{\infty}\{ (z_k)_{k \in N}  \in  R^{\infty} ~\&~ \underline{\lim}_{n \to \infty}T^{**}_n((z_k)_{1 \le k \le n})=\mu \}\\
& \ge \gamma_{\theta}^{\infty}\{ (z_k)_{k \in N}  \in  R^{\infty} ~\&~   \lim_{N \to \infty}T^{**}_n((z_k)_{1 \le k \le n})=\mu \}=1,\\
\end{align*}
which means  that  $T_{**}^{(0)}$ is an infinite-sample consistent  estimators of  the parameter $\theta$   for the family of probability measures$(\gamma_{\theta}^{\infty})_{\theta>0}$.

Similarly, we have

\begin{align*}
&\gamma_{\theta}^{\infty}\{ (z_k)_{k \in N} \in  R^{\infty} ~\&~  T_{**}^{(1)}((z_k)_{k \in N})=\mu \}\\
&=\gamma_{\theta}^{\infty}\{ (z_k)_{k \in N} \in  R^{\infty} ~\&~ \overline{\lim}_{n \to \infty}T^{**}_n((z_k)_{1 \le k \le n})=\mu\}\\
& \ge\gamma_{\theta}^{\infty}\{ (z_k)_{k \in N}  \in  R^{\infty} ~\&~   \lim_{N \to \infty}T^{**}_n((z_k)_{(z_k)_{1 \le k \le n}})=\mu \}=1,
\end{align*}
which means  that  $T_{**}^{(1)}$  is an infinite-sample consistent  estimators of  the parameter $\theta$   for the family of probability measures$(\gamma_{\theta}^{\infty})_{\theta>0}$.

\end{proof}

\noindent{\bf Remark 3.3.2}  Note that an
existence of  infinite sample consistent estimators  $T_{**}^{(0)}$  and  $T_{**}^{(1)}$  of  the rate  $\theta$   for the family of probability measures $(\gamma_{\theta}^{\infty})_{\theta>0}$
 (cf. Theorem 3.3.1) implies that the
family $(\gamma_{\theta}^{\infty})_{\theta>0}$   is strictly separated.

\subsection{Estimation of   the square of the  degree  of volatility  $\sigma$  around it caused by shocks  in  Ornstein - Uhlenbeck stochastic model }

We begin this subsection by the following proposition.

\noindent{\bf Lemma 3.4.1}  {\it  For  $t>0$, $x_0 \in R$, $\theta>0$, $\mu \in R$  and $\sigma>0$,   let's  $\gamma_{(t,x_0.\theta,\mu,\sigma)}$ be a Gaussian probability  measure in  $R$ with the mean  $m_t=x_0 e^{-\theta t}+\mu(1-e^{-\theta t}) $ and  the variance $\sigma_t^2=\frac{\sigma^2}{2\theta}\left( 1 - e^{-2\theta s}\right)$. Assuming that  parameters $x_0$, $t$, $\mu$  and $\theta$ are fixed.
For $\sigma^2 >0$,  let's  denote by $\gamma_{\sigma^2}$ the measure $\gamma_{(t,x_0,\theta,\mu,\sigma)}$.
Let define the estimate $T^{***}_n : R^n  \to R$    by the following formula
$$
T^{***}_n((z_k)_{1 \le k \le n})=\frac{2\theta \sum_{k=1}^n \left( z_k   -    x_0 e^{-\theta t}  -\mu (1-e^{-\theta t}) \right)^2}{n\left( 1 - e^{-2\theta s}\right)}.    \eqno  (3.4.1)
$$
Then we get
$$
\gamma_{\sigma^2}^{\infty}\{ (z_k)_{k \in N} :  (z_k)_{k \in N} \in  R^{\infty} ~\&~   \lim_{N \to \infty}T^{***}_n((z_k)_{1 \le k \le n})=\sigma^2 \}=1,\eqno  (3.4.2)
$$
provided that $T^{***}_n$ is a consistent estimator of  the square of the degree  of volatility  $\sigma$  around it caused by shocks  in the sense of convergence  almost everywhere   for the family of  probability  measures $(\gamma_{\sigma^2}^{\infty})_{\sigma^2 >0}$. }

\begin{proof} Let's consider probability space $(\Omega, \mathcal{F},P)$, where $\Omega=R^{\infty}$, $\mathcal{F}=B(R^{\infty})$, $P=\gamma_{\sigma^2}^{\infty}$.

For $k \in N$ we consider $k$-th projection $Pr_k$   defined on  $R^{\infty}$  by
$$
Pr_k( (x_i)_{i \in N})=x_k\eqno  (3.4.3)
$$
for $(x_i)_{i \in N} \in R^{\infty}$.

It is obvious that $(Pr_k)_{k \in N}$ is  sequence of independent Gaussian random variables  with the mean  $m_t=x_0 e^{-\theta t}+\mu(1-e^{-\theta t}) $ and the variance  $\sigma_t^2=\frac{\sigma^2}{2\theta}\left( 1 - e^{-2\theta s} \right)$. By use  Kolmogorov's  Strong  Law of Large numbers   for a sequence of independent  identically distributed  random variables  $(X_n)_{n \in N}$, where
$$
X_n( (z_j)_{j \in N})=\left (Pr_n( (z_j)_{j \in N}) -  x_0 e^{-\theta t}  -\mu (1-e^{-\theta t})\right)^2
$$
for $n \in N$, we get

$$\gamma_{\sigma^2}^{\infty}\{  (z_j)_{j \in N} \in  R^{\infty} ~\&~   \lim_{n \to \infty}\frac{\sum_{k=1}^n\left (Pr_k( (z_j)_{j \in N}) -  x_0 e^{-\theta t}  -\mu (1-e^{-\theta t})\right)^2 }{n}
$$
$$
=\frac{\sigma^2}{2\theta}\left( 1 - e^{-2\theta s}\right)  \}=1,  \eqno  (3.4.4)
$$
which implies
$$
\gamma_{\sigma^2}^{\infty} \{ (z_j)_{j \in N} \in  R^{\infty} ~\&~   \lim_{n \to \infty} \frac{2\theta \sum_{k=1}^n \left( z_k   -    x_0 e^{-\theta t}  -\mu (1-e^{-\theta t}) \right)^2}{n\left( 1 - e^{-2\theta s}\right)}=\sigma^2
  \}
$$
$$
=\gamma_{\sigma^2}^{\infty}\{ (z_j)_{j \in N}  \in  R^{\infty} ~\&~   \lim_{n \to \infty}T^{***}_n((z_j)_{1 \le j \le  n})=\sigma^2\}=1.
$$

\end{proof}

\noindent{\bf Remark 3.4.1}  By  use  Remark 2.1.1  and   Lemma 3.4.1  we deduce that  $T^{***}_n$ is a consistent estimator of the  parameter $\sigma^2$  in the sense of
convergence in probability for the family  of probability measures  $(\gamma_{\sigma^2}^{\infty})_{\sigma>0}$  as well $T^{***}_n$ is a consistent estimator  of the parameter $\sigma^2$   in the sense of
convergence in distribution for  the same family  of probability measures.

\medskip

\noindent{\bf Theorem 3.4.1}{\it ~ Suppose that the family of probability measures  $(\gamma_{\sigma^2}^{\infty})_{\sigma^2>0}$  and the estimators $T^{***}_n:R^n \to R (n \in N)$ come  from   Lemma 3.4.1.  Then the estimators
$T_{***}^{(0)}:R^{\infty} \to  R $  and  $T_{***}^{(1)}:R^{\infty} \to  R $ defined by
$$
T_{***}^{(0)}((z_k)_{k \in N})=\underline{\lim}_{n \to \infty}T^{**}_n((z_k)_{1\le k \le n}) \eqno(3.4.5)
$$
and
$$
T_{***}^{(1)}((z_k)_{k \in N})=\overline{\lim}_{n \to \infty}T^{**}_n((z_k)_{1\le k \le n}). \eqno(3.4.6)
$$
are  infinite-sample consistent  estimators  of the square of   the degree  of volatility  $\sigma$  around it caused by shocks   for the family of probability measures $(\gamma_{\sigma^2}^{\infty})_{\sigma^2>0}$ . }

\begin{proof} Note that  we have
\begin{align*}
&\gamma_{\sigma^2}^{\infty}\{ (z_k)_{k \in N}  \in  R^{\infty} ~\&~  T_{***}^{(0)}((z_k)_{k \in N})=\sigma^2 \}\\
&=\gamma_{\sigma^2}^{\infty}\{ (z_k)_{k \in N}  \in  R^{\infty} ~\&~ \underline{\lim}_{n \to \infty}T^{***}_n((z_k)_{1 \le k \le n})=\sigma^2 \}\\
& \ge \gamma_{\sigma^2}^{\infty}\{ (z_k)_{k \in N}  \in  R^{\infty} ~\&~   \lim_{N \to \infty}T^{***}_n((z_k)_{1 \le k \le n})=\sigma^2 \}=1,\\
\end{align*}
which means  that  $T_{***}^{(0)}$ is an infinite-sample consistent  estimators of  the parameter $\sigma^2$  for the family of probability measures  $(\gamma_{\sigma^2}^{\infty})_{\sigma^2>0}$.

Similarly, we have

\begin{align*}
&\gamma_{\sigma^2}^{\infty}\{ (z_k)_{k \in N} \in  R^{\infty} ~\&~  T_{***}^{(1)}((z_k)_{k \in N})=\sigma^2 \}\\
&=\gamma_{\sigma^2}^{\infty}\{ (z_k)_{k \in N} \in  R^{\infty} ~\&~ \overline{\lim}_{n \to \infty}T^{***}_n((z_k)_{1 \le k \le n})=\sigma^2\}\\
& \ge \gamma_{\sigma^2}^{\infty}\{ (z_k)_{k \in N}  \in  R^{\infty} ~\&~   \lim_{N \to \infty}T^{***}_n((z_k)_{(z_k)_{1 \le k \le n}})=\sigma^2\}=1,
\end{align*}
which means  that  $T_{***}^{(1)}$  is an infinite-sample consistent  estimators of the parameter $\sigma^2$    for the family of probability measures  $(\gamma_{\sigma^2}^{\infty})_{\sigma^2>0}$.

\end{proof}

\noindent{\bf Remark 3.4.2}  Note that an
existence of  infinite sample consistent estimators  $T_{***}^{(0)}$  and  $T_{***}^{(1)}$ of the parameter $\sigma^2$   for the family of probability measures
$(\gamma_{\sigma^2}^{\infty})_{\sigma^2>0}$
 (cf. Theorem 3.4.1) implies that the
family of probability measures   $(\gamma_{\sigma^2}^{\infty})_{\sigma^2>0}$  is strictly separated.

\section{Simulation  of  the Ornstein - Uhlenbeck stochastic  process   and  estimation  it's  parameters }

In this section  we  give a short explanation whether can be  obtained  the simulations  of the Ornstein - Uhlenbeck process. Similar simulations   can be found in  \cite{Smith2010}.

\begin{figure}[h]
\center{\includegraphics[width=0.8\linewidth]{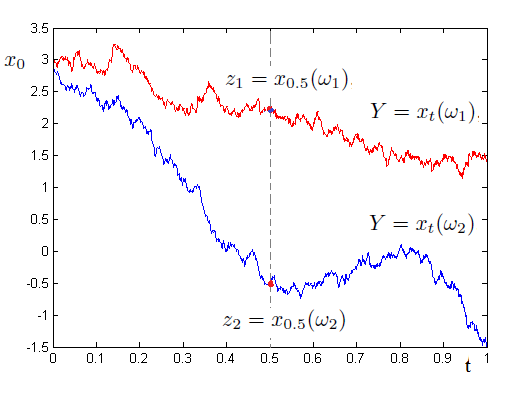}}
\caption{ Ornstein - Uhlenbeck's   two  trajectories   when  $\theta=0.5, ~\sigma=1,~\mu=-3$  and  $x_0=3.$
 }
\label{ris:image}
\end{figure}

The simulation  of  the Ornstein-Uhlenbeck process  can be  obtained  as follows:
$$ x_t=x_0 e^{-\theta t} +\mu (1-e^{-\theta t})+
{\sigma\over\sqrt{2\theta}}e^{-\theta t}W_{e^{2\theta t}-1}, \eqno  (4.1)$$
where $W_t$ denotes Wiener process.

Wiener (1923)  gave a representation of a Brownian path in terms of a random Fourier series. If  $(\xi_n)_{n \in N}$  is the sequence  of  independent standard Gaussian random  variables, then
$$W_t=\xi_0 t+ \sqrt{2}\sum_{n=1}^\infty\xi_n\frac{\sin \pi n t}{\pi n}  \eqno (4.2)$$
represents  a Brownian motion on $[0,1]$.

Following   Karhunen -Loeve well known  theorem  (see,  \cite{Loeve1948}, \cite{Karhunen1947} ),   the scaled process
$$\sqrt{c}\,  W\left(\frac{t}{c}\right) \eqno (4.3) $$
is a Brownian motion on $[0,c]$.

Let  $(y^{(k)}_n)_{n \in N}$ be a sequence of  real  numbers defined by $y^{(k)}_n=n \sqrt{p_k}-[ n \sqrt{p_k}]$ for $k,n \in N$, where $[ \cdot]$ denotes the integer part of the real number and $(p_n)_{n \in N}$ denotes the set of all prime numbers. Note that this sequence is uniformly distributed  in $(0,1)$  for each  $k \in N$ (see, for example, \cite{KuNi74}) .

Let $\Phi$ be a standard Gaussian  distribution function in $R$ . Then  following Lemma 2.1.1,  the  sequence  $(x^{(k)}_n)_{n \in N}=(\Phi^{-1}(y^{(k)}_n))_{n \in N}$ will be   $\gamma$-uniformly distributed  in  $R$  for each  $k \in N$, where $\gamma$ denotes a standard Gaussian measure in $R$.

In our simulation we  use  MatLab command  {\bf random('Normal', 0, 1, p, q)}   which  generates   $\gamma$ -uniformly  distributed    sequences  $(x^{(k)}_n)_{1 \le n \le q}  ~(1 \le k \le p)$.

\vspace{.08in} \noindent {\bf Table 4.1.  The value $z_k$  of the Ornstein-Uhlenbeck's   $k$-th  trajectory at moment $t=0.5$  when  $\theta=0.5, ~\sigma=1,~\mu=-3$  and  $x_0=3.$   }

~\begin{center}
\begin{tabular}{| l | l | l | l | l | l | l | l | l |  l |}

\hline $ k $&  $z_k$ &      $ k $ &  $z_k$ &        $ k $ & $z_k$ &      $ k $ & $z_k$&     $ k $ & $z_k$   \\

\hline

\hline  $1 $&  $2.7082$ &      $ 21 $ &  $1.2571$ &        $ 41 $ & $1.2185$ &      $61 $ & $1.9426$&     $ 81 $ & $2.7082$  \\

\hline

\hline  $2 $&  $2.0594$ &      $ 22 $ &  $2.1261$ &        $ 42 $ & $3.0131$ &      $62 $ & $1.453$&     $ 82 $ & $2.0594$  \\

 \hline

\hline  $3 $&  $2.1303$ &      $ 23 $ &  $2.6017$ &        $ 43 $ & $2.0324$ &      $63 $ & $1.6666$&     $ 83 $ & $2.1303$  \\

 \hline

\hline  $4 $&  $2.3939$ &      $ 24 $ &  $0.7975$ &        $ 44 $ & $1.8216$ &      $64 $ & $1.2806$&     $ 84 $ & $2.3939$  \\

 \hline

\hline  $5 $&  $2.641$ &      $ 25 $ &  $1.9225$ &        $ 45 $ & $1.1374$ &      $65 $ & $1.3268$&     $ 85 $ & $2.641$  \\

 \hline

\hline  $6 $&  $1.1519$ &      $ 26 $ &  $1.8187$ &        $ 46 $ & $2.7327$ &      $66 $ & $1.4312$&     $ 86 $ & $1.1519$  \\

 \hline

\hline  $7 $&  $1.6549$ &      $ 27 $ &  $1.8187$ &        $ 47 $ & $2.3649$ &      $67 $ & $2.7034$&     $ 87 $ & $1.6549$  \\

 \hline

\hline  $8 $&  $1.2017$ &      $ 28 $ &  $1.1202$ &        $ 48 $ & $1.3785$ &      $68 $ & $1.227$&     $ 88 $ & $0.6265$  \\

 \hline

\hline  $9 $&  $1.261$ &      $ 29 $ &  $0.3467$ &        $ 49 $ & $2.6211$ &      $69 $ & $1.0065$&     $ 89 $ & $1.2017$  \\

 \hline

\hline  $10 $&  $0.8576$ &      $ 30 $ &  $1.2734$ &        $ 50 $ & $1.258$ &      $70 $ & $0.7277$&     $ 90 $ & $1.261$  \\

 \hline

\hline  $11 $&  $1.3968$ &      $ 31 $ &  $2.6075$ &        $ 51 $ & $1.5606$ &      $71 $ & $1.3327$&     $ 91 $ & $0.8576$  \\

 \hline

\hline  $12 $&  $2.8304$ &      $ 32 $ &  $1.1872$ &        $ 52 $ & $2.0278$ &      $72 $ & $1.3528$&     $ 92 $ & $1.3968$  \\

 \hline

\hline  $13 $&  $1.1969$ &      $ 33 $ &  $1.9519$ &        $ 53 $ & $1.3095$ &      $73$ & $2.102$&     $ 93 $ & $2.8304$  \\

 \hline

\hline  $14 $&  $3.0469$ &      $ 34 $ &  $1.9615$ &        $ 54 $ & $1.8024$ &      $74 $ & $1.1705$&     $ 94 $ & $1.1969$  \\

 \hline

\hline  $15 $&  $0.7784$ &      $ 35 $ &  $1.6775$ &        $ 55 $ & $1.62$ &      $75 $ & $1.162$&     $ 95 $ & $3.0469$  \\

 \hline

\hline  $16 $&  $1.6111$ &      $ 36 $ &  $2.5195$ &        $ 56 $ & $0.9569$ &      $76 $ & $0.9056$&     $ 96 $ & $0.7784$  \\

 \hline

\hline  $17 $&  $1.1053$ &      $ 37 $ &  $1.894$ &        $ 57 $ & $0.8123$ &      $77$ & $0.6306$&     $ 97 $ & $1.6111$  \\

 \hline

\hline  $18 $&  $1.2695$ &      $ 38 $ &  $0.9174$ &        $ 58 $ & $0.9781$ &      $78$ & $0.3304$&     $ 98 $ & $1.1053$  \\

 \hline

\hline  $19 $&  $1.2756$ &      $ 39 $ &  $1.5291$ &        $ 59 $ & $1.9541$ &      $79$ & $1.0314$&     $ 99 $ & $1.2695$  \\

 \hline

\hline  $20 $&  $1.711$ &      $ 40 $ &  $1.3806$ &        $ 60 $ & $1.4921$ &      $80$ & $1.9173$&     $ 100 $ & $1.2756$  \\

\hline

\end{tabular}
\end{center}

Note that
$$
W^{(k)}_{e^{2\theta t}-1}=\xi^{(k)}_0 (e^{2\theta t}-1)+ \sqrt{2}\sum_{n=1}^\infty(\xi^{(k)}_n)\frac{\sin \pi n (e^{2\theta t}-1)}{\pi n} \eqno ( 4.4)
$$
will be the value of the  Wiener's  $k$-th  trajectory  at moment  $e^{2\theta t}-1$ for $k \in N$.

Hence the value of the  Ornstein-Uhlenbeck's   $k$-th  trajectory  at moment $t$ will be

$$z_k=x_0 e^{-\theta t} +\mu (1-e^{-\theta t})+
{\sigma\over\sqrt{2\theta}}e^{-\theta t} \left(  \xi^{(k)}_0  (e^{2\theta t}-1)+ \sqrt{2}\sum_{n=1}^\infty(\xi^{(k)}_n)\frac{\sin \pi n (e^{2\theta t}-1)}{\pi n} \right ) \eqno (4.5)$$
for each $k \in N$.

In our simulation we consider

$$z_k=x_0 e^{-\theta t} +\mu (1-e^{-\theta t})+
{\sigma\over\sqrt{2\theta}}e^{-\theta t} \left(     \xi^{(k)}_0  (e^{2\theta t}-1)+ \sqrt{2}\sum_{n=1}^{800}(\xi^{(k)}_n)\frac{\sin \pi n (e^{2\theta t}-1)}{\pi n}  \right ) \eqno (4.6)
$$

for  $1 \le k \le  100$.

Below  we present some numerical results  obtaining by using
MatLab  and  Microsoft Excel. In our  simulation

(i) ~$n$  denotes  the number of trials;

(ii) $x_0=3$  is  the underlying asset initial price;

(iii)   $\mu=-3$  is the equilibrium or mean value supported by fundamentals;

(iv)    $\sigma=1$  is   the degree of volatility around it caused by shocks;

(v)   $\theta=0.5$  is  the rate by which these shocks dissipate and the variable reverts towards the mean;

(vi)   $t=0.5$ is  the  moment of the observation on the Ornstein-Uhlenbeck's  trajectories;

(vii)          $z_k$  is  the value of the  Ornstein-Uhlenbeck's  $k$-th  trajectory at moment $t=0.5$(see,  Figure 1  and  Table 4.1).

\vspace{.08in} \noindent {\bf Table  4.2.  The  value of the statistic  $T_n$  for the  sample $(z_k)_{1 \le k \le n}(n=5i :1 \le i \le 20)$ from the Table 4.1   in  the  Ornstein-Uhlenbeck's  stochastic model   when  $\theta=0.5, ~\sigma=1,~\mu=-3$  and  $x_0=3.$           }

~\begin{center}
\begin{tabular}{|l|l|l|l|l|l|}

\hline $ n $       &    $T_n$   & $x_0$ &     $n$ &   $T_n$   & $x_0$\\

\hline

\hline $ 5 $       &    $3.916479949$   & $3$ &     $55$ &    $3.050042943$  & $3$\\

 \hline

\hline $ 10 $       &    $3.171013312$   & $3$ &     $60$ &   $2.999422576
$ & $3$\\

\hline

\hline $ 15 $       &    $3.189798604$   & $3$ &     $65$ &   $2.985749187
$   & $3$\\

\hline

\hline $ 20 $       &    $3.053011377$    & $3$ &     $70$ &   $2.963503799
$   & $3$\\

\hline

\hline $ 25 $       &    $3.059916865$  & $3$ &     $75$ &   $2.944638775
$   & $3$\\

\hline

\hline $ 30 $       &    $2.933186125$    & $3$ &     $80$ &   $2.891140712
$  & $3$\\

\hline

\hline $ 35 $       &    $2.98020897$  & $3$ &     $85$ &   $2.951454785
$   & $3$\\

\hline

\hline $ 40 $       &    $2.978720876$   & $3$ &     $90$ &   $2.918940576
$   & $3$\\

 \hline

\hline $ 45 $       &    $3.005595171$   & $3$ &     $95$ &   $2.936244133
$   & $3$\\

\hline

\hline $ 50 $       &    $3.056170079$   & $3$ &   $100$ &    $2.90958959
$  & $3$\\

\hline

\end{tabular}
\end{center}

\medskip

\noindent{\bf Remark 4.1}  By  use  results  of  calculations   placed in the Table 4.2,  we see that the  consistent estimator $T_n$  works   successfully.

\medskip

\vspace{.08in} \noindent {\bf Table 4.3.  The  value of the statistic  $T^*_n$  for the  sample $(z_k)_{1 \le k \le n}(n=5i :1 \le i \le 20)$  from the Table 4.1    in  the  Ornstein-Uhlenbeck's  stochastic model   when  $\theta=0.5, ~\sigma=1,~\mu=-3$  and  $x_0=3.$   }
~\begin{center}
\begin{tabular}{|l|l|l|l|l|l|}

\hline $ n $       &    $T^*_n$   & $\mu$ &     $n$ &   $T^*_n$   & $\mu$\\

\hline

\hline $ 5 $       &    $0.226753293$   & $-3$ &     $55$ &    $-2.823808222$  & $-3$\\

 \hline

\hline $ 10 $       &    $-2.397894335$   & $-3$ &     $60$ &   $-3.002033002$ & $-3$\\

\hline

\hline $ 15 $       &    $-2.331754861$   & $-3$ &     $65$ &   $-3.05017443$   & $-3$\\

\hline

\hline $ 20 $       &    $-2.813356927$    & $-3$ &     $70$ &   $-3.12849625$   & $-3$\\

\hline

\hline $ 25 $       &    $-2.789044002$  & $-3$ &     $75$ &   $-3.194916445$   & $-3$\\

\hline

\hline $ 30 $       &    $-3.235239072$    & $-3$ &     $80$ &   $-3.38327305$  & $-3$\\

\hline

\hline $ 35 $       &    $-3.06968049$  & $-3$ &     $85$ &   $-3.170918559$   & $-3$\\

\hline

\hline $ 40 $       &    $-3.074919788$   & $-3$ &     $90$ &   $-3.285394965$   & $-3$\\

 \hline

\hline $ 45 $       &    $-2.980300456$   & $-3$ &     $95$ &   $-3.224472402$   & $-3$\\

\hline

\hline $ 50 $       &    $-2.80223573$   & $-3$ &   $100$ &    $-3.318318028$  & $-3$\\

\hline

\end{tabular}
\end{center}

\medskip

\noindent{\bf Remark 4.2}  By  use  results  of  calculations   placed in the Table 4.3,  we see that the  consistent estimator $T^*_n$  works   successfully.

\vspace{.08in} \noindent {\bf Table 4.4.  The  value of the statistic  $T^{**}_n$  for the  sample$(z_k)_{1 \le k \le n}(n=5i :1 \le i \le 20)$ from the Table 4.1  in  the  Ornstein-Uhlenbeck's  stochastic model   when  $\theta=0.5, ~\sigma=1,~\mu=-3$  and  $x_0=3.$   }

~\begin{center}
\begin{tabular}{|l|l|l|l|l|l|}

\hline $ n $       &    $T^{**}_n$   & $\theta$ &     $n$ &   $T^{**}_n$   & $\theta$\\

\hline

\hline $ 5 $       &    $0.215705016$   & $0.5$ &     $55$ &    $0.483388203$  & $0.5$\\

 \hline

\hline $ 10 $       &    $0.44379283$   & $0.5$ &     $60$ &   $0.500192489$ & $0.5$\\

\hline

\hline $ 15 $       &    $0.437713842$   & $0.5$ &     $65$ &   $0.504755927$   & $0.5$\\

\hline

\hline $ 20 $       &    $0.482407151$    & $0.5$ &     $70$ &   $0.512202556$   & $0.5$\\

\hline

\hline $ 25 $       &    $0.480126781$  & $0.5$ &     $75$ &   $0.518539409$   & $0.5$\\

\hline

\hline $ 30 $       &    $0.522396228$    & $0.5$ &     $80$ &   $0.536619648$  & $0.5$\\

\hline

\hline $ 35 $       &    $0.50660792$  & $0.5$ &     $85$ &   $0.516247561$   & $0.5$\\

\hline

\hline $ 40 $       &    $0.507105654$   & $0.5$ &     $90$ &   $0.527203992$   & $0.5$\\

 \hline

\hline $ 45 $       &    $0.498135817$   & $0.5$ &     $95$ &   $0.521365679$   & $0.5$\\

\hline

\hline $ 50 $       &    $0.481363743$   & $0.5$ &   $100$ &    $0.530366173$  & $0.5$\\

\hline

\end{tabular}
\end{center}

\medskip

\noindent{\bf Remark 4.3}  By  use  results  of  calculations   placed in the Table 4.4,  we see that the  consistent estimator $T^{**}_n$  works   successfully.

\vspace{.08in} \noindent {\bf Table 4.5.  The  value of the statistic  $T^{***}_n$  for the  sample $(z_k)_{1 \le k \le n}(n=5i :1 \le i \le 20)$ from the Table 4.1    in  the  Ornstein-Uhlenbeck's  stochastic model   when  $\theta=0.5, ~\sigma=1,~\mu=-3$  and  $x_0=3.$   }

~\begin{center}
\begin{tabular}{|l|l|l|l|l|l|}

\hline $ n $       &    $T^{***}_n$   & $\sigma^2$ &     $n$ &   $T^{***}_n$   & $\sigma^2$\\

\hline

\hline $ 5 $       &    $1.468059434$   & $1$ &     $55$ &    $1.019738082$  & $1$\\

 \hline

\hline $ 10 $       &    $1.071475$   & $1$ &     $60$ &   $1.013011822$ & $1$\\

\hline

\hline $ 15 $       &    $1.448094876$   & $1$ &     $65$ &   $0.950520018$   & $1$\\

\hline

\hline $ 20 $       &    $1.168384329$    & $1$ &     $70$ &   $0.979072751$   & $1$\\

\hline

\hline $ 25 $       &    $1.145106531$  & $1$ &     $75$ &   $0.944824889$   & $1$\\

\hline

\hline $ 30 $       &    $1.174454042$    & $1$ &     $80$ &   $1.011196765$  & $1$\\

\hline

\hline $ 35 $       &    $1.098947717$  & $1$ &     $85$ &   $1.03807104$   & $1$\\

\hline

\hline $ 40 $       &    $1.053231322$   & $1$ &     $90$ &   $1.030042542$   & $1$\\

 \hline

\hline $ 45 $       &    $1.074061413$   & $1$ &     $95$ &   $1.088067168$   & $1$\\

\hline

\hline $ 50 $       &    $1.106961814$   & $1$ &   $100$ &    $1.070420297$  & $1$\\

\hline

\end{tabular}
\end{center}

\medskip

\noindent{\bf Remark 4.3}  By  use  results  of  calculations   placed in the Table 4.5,  we see that the  consistent estimator $T^{***}_n$  works   successfully.

\end{document}